\documentclass{amsart}
\usepackage{ifthen, amssymb, tikz, natbib}
\usepackage[multiple]{footmisc}

\newboolean{oneResultSequence}
\setboolean{oneResultSequence}{false} 

\newcommand{\numberResults}{
\ifthenelse{\boolean{oneResultSequence}}{
\newtheorem{axiom}{Axiom}
\newtheorem{theorem}{Theorem}
\newtheorem{lemma}[theorem]{Lemma}
\newtheorem{proposition}[theorem]{Proposition}
\newtheorem{corollary}[theorem]{Corollary}
\theoremstyle{definition}
\newtheorem{definition}[theorem]{Definition}
\newtheorem{example}[theorem]{Example}
\newtheorem{xca}[theorem]{Exercise}
\theoremstyle{remark}
\newtheorem{remark}[theorem]{Remark}
\newtheorem{conjecture}[theorem]{Conjecture}
}{
\newtheorem{axiom}{Axiom}
\newtheorem{theorem}{Theorem}
\newtheorem{lemma}{Lemma}
\newtheorem{proposition}{Proposition}
\newtheorem{corollary}{Corollary}

\theoremstyle{definition}
\newtheorem{definition}{Definition}
\newtheorem{example}{Example}

\theoremstyle{remark}

}}

\def\lema#1#2{\begin{lemma} #2 \label{lem:#1} \end{lemma}}

\def\corol#1#2{\begin{corollary} #2 \label{cor:#1} \end{corollary}}

\def\propo#1#2{\begin{proposition} #2 \label{prp:#1} \end{proposition}}

\def\display#1#2{\begin{equation} #2 \label{eqn:#1} \end{equation}}
\def\eqn#1{(\ref{eqn:#1})}

\newboolean{RANW}
\setboolean{RANW}{false}
\newlength{\testWidth} \newlength{\badWidth}

\newcommand{\useRsltAndNumberwithin}{
\setboolean{RANW}{true}
\settowidth\badWidth{\ref{useRsltAndNumberwihin_generates_spurious_undefined-reference_messages.}}
\newcommand{\badMatch}{\equal{\the\testWidth}{\the\badWidth}}
}

\newcommand{\rslt}[1]{%
\ifthenelse{\boolean{RANW}}%
{%
{\settowidth{\testWidth}{\ref{lem:#1}}\ifthenelse{\badMatch}%
{%
\settowidth{\testWidth}{\ref{prp:#1}}\ifthenelse{\badMatch}%
{%
\settowidth{\testWidth}{\ref{thm:#1}}\ifthenelse{\badMatch}%
{%
\settowidth{\testWidth}{\ref{cor:#1}}\ifthenelse{\badMatch}{\textbf{result ??}}%
{corollary \ref{cor:#1}}%
}%
{theorem \ref{thm:#1}}%
}%
{proposition \ref{prp:#1}}%
}%
{lemma \ref{lem:#1}}%
}}%
{%
\ifthenelse{\ref{lem:#1} > 0}{lemma \ref{lem:#1}}%
{%
\ifthenelse{\ref{prp:#1} > 0}{proposition \ref{prp:#1}}%
{%
\ifthenelse{\ref{thm:#1} > 0}{theorem \ref{thm:#1}}%
{%
\ifthenelse{\ref{cor:#1} > 0}{corollary \ref{cor:#1}}{\textbf{Result ??}}%
}%
}%
}%
}}

\newcommand{\Rslt}[1]{%
\ifthenelse{\boolean{RANW}}%
{%
{\settowidth{\testWidth}{\ref{lem:#1}}\ifthenelse{\badMatch}%
{%
\settowidth{\testWidth}{\ref{prp:#1}}\ifthenelse{\badMatch}%
{%
\settowidth{\testWidth}{\ref{thm:#1}}\ifthenelse{\badMatch}%
{%
\settowidth{\testWidth}{\ref{cor:#1}}\ifthenelse{\badMatch}{\textbf{Result ??}}%
{Corollary \ref{cor:#1}}%
}%
{Theorem \ref{thm:#1}}%
}%
{Proposition \ref{prp:#1}}%
}%
{Lemma \ref{lem:#1}}%
}}%
{%
\ifthenelse{\ref{lem:#1} > 0}{Lemma \ref{lem:#1}}%
{%
\ifthenelse{\ref{prp:#1} > 0}{Proposition \ref{prp:#1}}%
{%
\ifthenelse{\ref{thm:#1} > 0}{Theorem \ref{thm:#1}}%
{%
\ifthenelse{\ref{cor:#1} > 0}{Corollary \ref{cor:#1}}{\textbf{Result ??}}%
}%
}%
}%
}}

\newcommand{\rsltref}[1]{%
\ifthenelse{\boolean{RANW}}%
{%
{\settowidth{\testWidth}{\ref{lem:#1}}\ifthenelse{\badMatch}%
{%
\settowidth{\testWidth}{\ref{prp:#1}}\ifthenelse{\badMatch}%
{%
\settowidth{\testWidth}{\ref{thm:#1}}\ifthenelse{\badMatch}%
{%
\settowidth{\testWidth}{\ref{cor:#1}}\ifthenelse{\badMatch}{\textbf{result ??}}%
{\ref{cor:#1}}%
}%
{\ref{thm:#1}}%
}%
{\ref{prp:#1}}%
}%
{\ref{lem:#1}}%
}}%
{%
\ifthenelse{\ref{lem:#1} > 0}{\ref{lem:#1}}%
{%
\ifthenelse{\ref{prp:#1} > 0}{\ref{prp:#1}}%
{%
\ifthenelse{\ref{thm:#1} > 0}{\ref{thm:#1}}%
{%
\ifthenelse{\ref{cor:#1} > 0}{\ref{cor:#1}}{\textbf{Result ??}}%
}%
}%
}%
}}

\def\Section#1#2{\section{#2\label{sec:#1}}}

\def\sec#1{{s}ection \ref{sec:#1}} 

\newenvironment{enum} 
{\begin{list}{\makebox[\labelwidth][l]{(\arabic{enumi})}}{\usecounter{enumi}}
\setcounter{enumi}{\value{equation}}}
{\setcounter{equation}{\value{enumi}} \end{list}}

\newcommand{\meti}[2]{\item #2 \label{eqn:#1}} 

\newcommand{\abbrevEnvir}{
\expandafter\newcommand\expandafter{\csname bi\endcsname}{\begin{itemize}} 
\expandafter\newcommand\expandafter{\csname ei\endcsname}{\end{itemize}}
\expandafter\newcommand\expandafter{\csname be\endcsname}{\begin{enumerate}} 
\expandafter\newcommand\expandafter{\csname ee\endcsname}{\end{enumerate}}
\expandafter\newcommand\expandafter{\csname bc\endcsname}{\begin{center}} 
\expandafter\newcommand\expandafter{\csname ec\endcsname}{\end{center}}
}

\newcommand{\from}{\colon}
\newcommand{\Nat}{\mathbb{N}}

\newcommand{\authorGreen}{\ifthenelse{\isundefined{\authorInPageHeader}}
{\author[\relax]{Edward J. Green}}
{\author[\authorInPageHeader]{Edward J. Green}}}

\authorGreen
\address{Department of Economics, The Pennsylvania State University,
  University Park, PA 16802, USA}
\email{eug2@psu.edu}

\newcommand{\num}{\varepsilon} 
\newcommand{\numtwo}{\theta} 
\newcommand{\cod}{\gamma} 
\newcommand{\codtwo}{\zeta}
\newcommand{\cpair}{\delta} 
\newcommand{\cpairtwo}{\theta} 
\newcommand{\merge}{\kappa} 

\newcommand{\field}{\mathcal{F}}

\newcommand{\Z}{\mathbb{Z}}

\newcommand{\id}{I}
\newcommand{\rcod}[1][\cod]{R_{#1}}
\newcommand{\scod}[1][\cod]{S_{#1}}
\newcommand{\tcod}[1][\cod]{T_{#1}}

\newcommand{\eqva}{F}
\newcommand{\eqvb}{G}
\newcommand{\eqvc}{H}

\newcommand{\proj}[2][2]{\pi^{#1}_{#2}} 
\newcommand{\Proj}[1]{\overline{\pi}_{#1}}
\newcommand{\seq}[1][{}]{\sigma^{#1}} 
\newcommand{\bound}{\beta}

\newcommand{\inv}[1]{#1^{-1}}

\newcommand{\tc}[1]{#1^+} 

\newcommand{\relpow}[2]{#1^{(#2)}} 

\newcommand{\head}{\eta}
\newcommand{\tail}{\tau}

\newcommand{\mand}{\text{\ and\ }}

\newcommand{\coarsep}{IC} \newcommand{\coarse}{\coarsep\ }
\newcommand{\finep}{FC} \newcommand{\fine}{\finep\ }

\providecommand{\citet}[1]{#1}
\providecommand{\citep}[1]{(\citet{#1})}

\setboolean{oneResultSequence}{false}
\numberResults 
\abbrevEnvir

\begin{document}


\title[semi-decidable equivalence relations]{Semi-decidable equivalence relations\\
obtained by composition and lattice join\\ 
of decidable equivalence relations}

\date{2018-01-17}

\begin{abstract}
Composition and lattice join (transitive closure of a union) of
equivalence relations are operations taking pairs of decidable
equivalence relations to relations that are semi-decidable, but not
necessarily decidable. This article addresses the question, is
\emph{every} semi-decidable equivalence relation obtainable in those
ways from a pair of decidable equivalence relations? It is shown that
every semi-decidable equivalence relation, of which every equivalence
class is infinite, is obtainable as both a composition and a lattice
join of decidable equivalence relations having infinite equivalence
classes. An example is constructed of a semi-decidable, but not
decidable, equivalence relation having finite equivalence classes that
can be obtained from decidable equivalence relations, both by
composition and also by lattice join. Another example is constructed,
in which such a relation cannot be obtained from decidable equivalence
relations in either of the two ways.
\end{abstract}

\maketitle

\Section{1}{Introduction}

Pullback, intersection, composition and lattice join of equivalence
relations are the operations by which, in practice, new equivalence
relations are typically formed from antecedent ones.\footnote{The
  transitive closure of a union of equivalence relations is the join
  of those relations (and their intersection is their meet) in the
  lattice of equivalence relations, in which a refinement of a
  relation is ordered below it.} All four operations preserve
semi-decidability.\footnote{In this article, \emph{computable} and
  \emph{decidable} describe the functions and relations, respectively,
  that have been called \emph{recursive} in the older terminology
  adopted by \citet{Kleene-1952} and \citet{Rogers-1967}.
  \emph{Partially computable} and \emph{semi-decidable} describe
  objects that Rogers called \emph{partially recursive} and
  \emph{recursively enumerable}.  Some authors also use the adjective,
  \emph{positive,} and the acronym, \emph{ceer,} to refer to a
  semi-decidable (or computably enumerable) equivalence
  relation. Concepts and results that will be introduced below without
  definition or proof can be found (in identical or transparently
  equivalent form) in the early chapters of Rogers' book. Although
  their article focuses on the specific topic of universal relations,
  \citet{AndrewsEtAl-2017}, broadly cover the research literature on
  semi-decidable equivalence relations in their
  bibliography.}\footnote{In this assertion and the following one, a
  pullback by a computable function is assumed.} However, while
pullback and intersection of decidable relations yield relations that
are also decidable, composition and lattice join take pairs of
decidable equivalence relations to relations that are not necessarily
decidable.  Since those two operations are defined by existential
quantification (over elements of their arguments, for composition; and
over sequences of those elements, for lattice join), an analogy with
Kleene's projection theorem---that a set is semi-decidable iff it is
defined by existential quantification over a decidable set---suggests
the possibility that \emph{every} semi-decidable equivalence relation
might be obtainable as a composition or lattice join of decidable
equivalence relations. This article investigates that conjecture. It
is shown in \rslt{18} that every semi-decidable equivalence relation,
of which every equivalence class is infinite, is so obtainable from
decidable equivalence relations having infinite equivalence
classes. \Rslt{23} specifies a semi-decidable, but not decidable,
equivalence relation, having finite equivalence classes, that can be
obtained from decidable equivalence relations in each of those two
ways. An example is constructed also, in \rslt{25} of \rslt{24}, of
such a relation that cannot be obtained from decidable equivalence
relations in either way.

\Section{2}{Computability and decidability, equivalence
relations}

Let $\Nat$ denote $\{ 0,1,2\dots \}$, let $\Nat_+$ denote $\{
1,2,3\dots \}$, and let $\Z$ denote the integers. In quantified
formulae below, variables will range over $\Nat$. \emph{Enumeration of
  $R$} will be used in this article to mean \emph{computable function
  from $\Nat_+$ onto $R$.}\footnote{The choice to make $\Nat_+$ the
  domain, rather than $\Nat$, is motivated by the definition of a walk
  later in this section, where it will be convenient that $n \neq -n$
  for all $n$ in the domain of the enumeration.} Recall that $R$ is
semi-decidable iff there is an enumeration of $R$.

For every $1 \le i \le n \in \Nat_+$, there is a computable function
$\proj[n]{i} \from \Nat^n \to \Nat$ such that, for every $e =
(x_1,\dotsc,x_n) \in \Nat^n$,
\display{1}{\proj[n]{i}(e) = x_i}
Due to this fact, there is no loss of generality in restricting
attention to relations that are subsets of $\Nat \times \Nat$ rather
than studying subsets of $\Nat^n \times \Nat^n$ in general.
(Cf.\ \citet[pp.\ 64--66]{Rogers-1967}.) This simplification will be
made henceforth in this article, recognizing that all results proved
here can be generalized.

An \emph{equivalence relation} on $\Nat$ is a reflexive, transitive,
symmetric relation, represented as a subset of $\Nat \times \Nat$.
Define an equivalence relation, $E$, to be \emph{\coarse} if every
equivalence class of $E$ has the cardinality of $\Nat$, and to be
\emph{\fine} if every equivalence class is finite.\footnote{`\finep',
  acronym for `finite class', was introduced by
  \citet{GaoGerdes-2001}. Correspondingly, `\coarsep' stands for
  `infinite class'.} Let $[i]_E$ (or simply $[i]$, if the meaning is
clear) denote the equivalence class of $i$ in $E$. The \emph{partition
  corresponding to} $E$ is the set of those equivalence classes.

Define the \emph{field} of a binary relation, $H$, by\footnote{For
  $(i,j) \in \Nat \times \Nat$, define $\inv{(i,j)} = (j,i)$. For $R
  \subset \Nat \times \Nat$, define $\inv R = \{ \inv e \mid e \in R
  \}$, the \emph{converse relation of $R$}.}
\display{2}{\field(H) = \{ i \mid\; \exists j \, (i,j) \in H \cup \inv
  H \} }
The notion of an equivalence class can be extended to symmetric,
transitive relations (with $[i]_H = \emptyset$ if $i \notin
\field(H)$) by defining
\display{3}{[i]_H = \{ j \mid (i,j) \in H \} }
A symmetric, transitive relation, $H$, will be said to be \coarse if
$[i]$ is infinite for every $i$ in $\field(H)$. The following, obvious
lemma will be used in \sec{7}.

\lema{1}{If $H$ is a symmetric, transitive relation on
  $\Nat$, then $i \in \field(H)$ iff $(i,i) \in H$.
If each of $H$ and $J$ is a symmetric, transitive relation on
  $\Nat$ and is \coarsep, and if $\field(J) = \Nat \setminus
  \field(H)$, then $H \cup J$ is an \coarse equivalence relation on
  $\Nat$.}

Proofs below will require a formal definition of the transitive
closure of a binary relation and the statement of an equivalent
characterization of it (\rslt{3}) for symmetric, semi-decidable relations. Let
$\tc R$ denote the transitive closure of $R$. Define $\relpow{R}{1} =
R$ and $\relpow{R}{n+1} = R \relpow{R}{n}$.\footnote{$RS$ denotes the
  composition of $R$ and $S$.} Then
\display{4}{\tc R = \bigcup_{n \in \Nat_+} \relpow R n}
The following lemma, proved by showing inductively that the hypothesis
entails that $\relpow{R}{n} \subseteq R$, will be used in \sec{8}.

\lema{2}{Suppose that the equivalence class in $E$ of every number is
  a singleton or a pair. If $R \subseteq E$ is reflexive, then $\tc R
  = R$. If $S \subseteq E$ is reflexive and symmetric, then $S$ is an
  equivalence relation.}

Suppose that $R \subset \Nat \times \Nat$ is a semi-decidable relation
enumerated by $\num$. For $k \in \Nat_+$, define
\display{5}{\tail(k) = \head(-k) = \proj 1 \num(k) \mand \head(k) =
  \tail(-k) = \proj 2 \num(k)}
That is, the ordered pairs in $R \cup \inv R$ are viewed as edges of a
directed graph, $\num$ is extended to $\Z \setminus \{ 0 \}$ by defining $\num(-k) = \inv{\num(k)}$, and $x
\in \Z \setminus \{ 0 \}$ is interpreted as being a directed edge with
tail $\tail(x)$ and head $\head(x)$.

A \emph{walk} (of length $n$) from $i$
to $j$ in $R$ is a sequence, $(x_1,\dotsc,x_n) \in (\Z \setminus \{ 0
\})^n$ such that $\tail(x_1) = i$, $\head(x_n) = j$, and, for $1 \le k
< n$, $\head(x_k) = \tail(k+1)$.

The next lemma follows straightforwardly by induction from these
definitions.

\lema{3}{If $\num$ enumerates $R$, then $(i,j) \in \tc{(R \cup \inv R)}$
  iff there is a walk from $i$ to $j$ in $R$. If there is a walk of
  length $n$, then $(i,j) \in \relpow{(R \cup \inv R)}{n}$.}

\Section{3}{Coding a semi-decidable equivalence relation}

Henceforth throughout this article, it will be assumed that
\begin{enum}
\meti{6}{$E \subseteq \Nat^2$ is an equivalence relation and $\num
  \from \Nat_+ \to E$ enumerates $E$.}
\end{enum}

In this section, a computable injection, $\cod \from \Nat_+ \to
\Nat_+$, will be defined that will be seen to encode $E$ implicitly,
where $E$ is a semi-decidable, \coarse equivalence relation.  Defining
\display{7}{\cpair = \num \cod,}
Later in the article, $\cpair$ will be shown to enumerate a decidable
subset of $E$. The transitive closure of $\cpair(\Nat) \cup
\inv{(\cpair(\Nat))}$, to be denoted by $\eqva$, will be shown to be a
decidable, \coarse equivalence relation. $E$ will subsequently be
characterized as $\eqva\eqvb\eqva$, where $\eqvb$ is another such
relation that is also derived from $\cod$. Because $\eqva \subset E$
and $\eqvb \subseteq E$, it follows that $E$ is the join of $\eqva$
and $\eqvb$ in the lattice of equivalence relations.

A computable function, $\cod \from \Nat_+ \to \Nat_+$, satisfying
conditions \eqn{8} below will be called a \emph{coding of $\num$,} and
will be said to \emph{code} $\num$. It will be called a \emph{coding
  of $E$} iff it codes some enumeration of $E$. An equivalence
relation, of which some enumeration has a coding, will be called
\emph{codable}. \Rslt{8} will establish that a semi-decidable
equivalence relation is codable iff \emph{every} enumeration has a
coding, and also iff the relation is \coarsep.
\display{8}{\begin{gathered}
    1 < \cod(1) \qquad \proj 1 \cpair(1) = \proj 1 \num(1) \qquad
    \max \{ 1, \proj 1 \num(1) \} < \proj 2 \cpair(1)\\
    \cod(n) < \cod(n+1) \qquad \proj 1 \num(\cod(n+1)) =
    \proj 1 \num(n+1)\\
    \max \{ \proj 1 \cpair(n+1), \proj 2 \cpair(n) \} < \proj 2
    \cpair(n+1)) \}
\end{gathered}}

\lema{4}{If $\cod \from \Nat_+ \to \Nat_+$, satisfies conditions \eqn{8},
  and if equation \eqn{7} defines $\cpair$, then $\cod$ and $\cpair$ are
  computable.}

\begin{proof}
Conditions \eqn{8} define a partially computable function, $\cod$, by
recursion. By hypothesis, $\cod$ is total. A partially computable,
total function is computable, so $\cod$ is computable. A composition
of computable functions is computable, so $\cpair$ is computable.
\end{proof}

The relations specified in \eqn{8} are depicted in the following
diagram, figure \eqn{9}. The higher endpoint of a solid line segment is explicitly
specified to be a larger number than the lower endpoint
is. Additionally, from $1 < \cod(1)$ and $\cod(n) < \cod(n+1)$, it
follows by induction that $n < \cod(n)$ and $n+1 < \cod(n+1)$. By
parallel reasoning, $n < \proj 2 \cpair(n)$ and $n+1 < \proj 2
\cpair(n+1)$.\footnote{These parallel conclusions
regarding $\cod$ and $\cpair$ will be restated in
    \rslt{5}.} Also $n < n+1$. The heights of the opposite endpoints
  of the respective dashed line segments indicate these relative
  magnitudes. A dotted line segment indicates that its endpoints are
  related as argument and image of a function (which may be a
  composite function). The relative magnitudes of opposite endpoints
  of dotted line segments are not constrained by \eqn{8}.

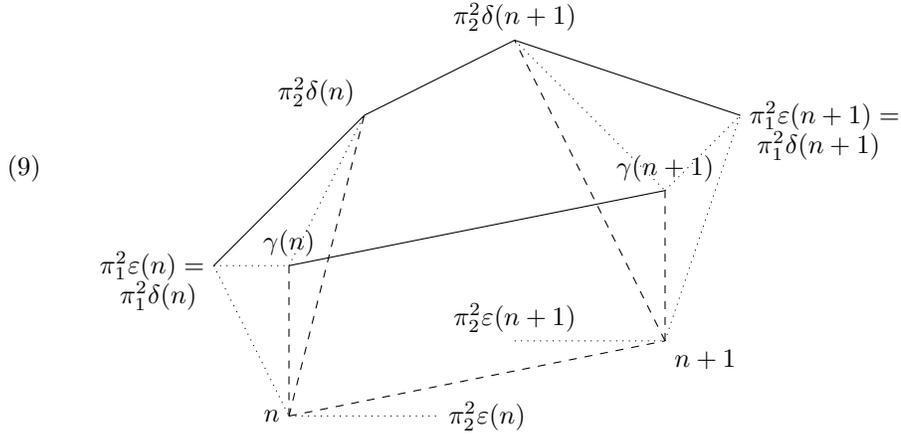
\strut

\display{9}{
\begin{tikzpicture}[baseline=(a.base)]
\draw[dotted] (2,3) node[left]{$\proj 1 \num(n) =$} -- (3,1)
node[left]{$n$};
\draw (1.9,2.6) node[left]{$\proj 1 \cpair(n)$};
\draw[dotted] (2,3) -- (3,3) node[above]{$\cod(n)$};
\draw[dotted] (3,1) -- (5,1) node[right]{$\proj 2 \num(n)$};
\draw[dotted] (3,3) -- (4,5) node[above left]{$\proj 2 \cpair(n)$};
\draw (2,3) -- (4,5);
\draw[dashed] (3,1) -- (3,3);
\draw[dashed] (3,1) -- (4,5);
\draw (4,5) -- (6,6) node[above]{$\proj 2 \cpair(n+1)$};
\draw[dashed] (3,1) -- (8,2) node[below right]{$n+1$};
\draw[dashed] (8,2) -- (6,6);
\draw[dotted] (6,2) node[above]{$\proj 2 \num(n+1)$} -- (8,2);
\draw[dashed] (8,2) -- (8,4) node[above](a){$\cod(n+1)$};
\draw[dotted] (8,2) -- (9,5) node[right]{$\proj 1 \num(n+1) =$};
\draw (9.1,4.6) node[right]{$\proj 1 \cpair(n+1)$};
\draw (6,6) -- (9,5);
\draw[dotted] (6,6) -- (8,4);
\draw (3,3) -- (8,4);
\draw[dotted] (8,4) -- (9,5);
\end{tikzpicture}
}

\strut

\noindent Two crucial features of this diagram encapsulate the role
that the construction of $\cod$ will play in the following proofs:
\begin{align*}
\text{The relative\ }&\text{magnitudes of\ } \proj 1 \num(n) \mand \proj 2
\num(n) \text{\ are indeterminate,}\\ \text{but\ } &
\proj 1 \cpair(n) < \proj 2 \cpair(n)\\
\text{The relative\ }&\text{magnitudes of\ } \proj 2 \num(n) \mand \proj 2
\num(n+1) \text{\ are indeterminate,}\\ \text{but\ } &
\proj 2 \cpair(n) < \proj 2 \cpair(n+1)
\end{align*}

\lema{5}{If $\cod$ codes $\num$, then both $\cod$ and $\proj 2
  \cpair$ are strictly increasing functions. For all $n$, $n <
  \cod(n)$ and $n < \proj 2 \cpair(n)$.}

\begin{proof}

Equation \eqn{8} requires explicitly that $\cod$ must be increasing,
and also that $1 < \cod(1)$. If $n < \cod(n)$, then $n+1 < \cod(n) + 1
\le \cod(n+1)$, so $\cod$ satisfies $\forall n \; n < \cod(n)$ by
induction. Taking $k = \cod(1)$ in the first equation and $k =
\cod(n+1)$ in the second equation of \eqn{8}, those equations specify
that $1 < \proj 2 \num \cod(1) = \proj 2 \cpair(1)$ and that $\proj 2
\cpair(n) < \proj 2 \num \cod(n+1) = \proj 2 \cpair(n+1)$. The second
of these inequalities states that $\proj 2 \cpair$ is strictly
increasing, and the two inequalities together imply by induction that
$\forall n \; n < \proj 2 \cpair(n)$.
\end{proof}

\lema{6}{If some enumeration, $\num$, of an equivalence relation, $E$,
  has a coding, $\cod$, then, for all $i \in \Nat_+$, $\{ n \mid \proj 1
  \cpair(n) = i \}$ is infinite, and $E$ is \coarsep. Specifically,
  where $\num(k) = (i,i)$, $\{\cpair\cod^n(k) \mid n \in \Nat \} =
  \{ i \} \times \{\proj 2 \cpair\cod^n(k) \mid n \in \Nat \}$ is
  an infinite subset of $E$.}

\begin{proof}
Assume that $\cod$ codes $\num$, an enumeration of $E$, and that
$\cpair = \num\cod$. Let $i \in \Nat_+$, and let
$\num(k) = (i,i)$. By induction, for all $n$, $(\proj 1 \num
\cod^n(k), \proj 2 \num \cod^n(k)) = (i, \proj 2 \num \cod^n(k))
\in E$. By induction, invoking \rslt{5}, for all $n$, $\cod^n(k) <
\cod(\cod^n(k)) = \cod^{n+1}(k)$ and consequently $\proj 2
\cpair(\cod^n(k)) < \proj 2 \cpair(\cod^{n+1}(k))$.\footnote{For any
  mapping $f \from \Nat_+ \to \Nat_+$, $f^0$ denotes the identity
  mapping and $f^{n+1} = ff^n$.}

Thus, for every $i$, $\{ \proj 2 \cpair\cod^n(k) \mid n \in \Nat \}$ is an
infinite subset of $[i]$, where $\num(k) = (i,i)$.
\end{proof}

\lema{7}{If $E$ is an \coarse equivalence relation with
  enumeration $\num$, then there exists a coding, $\cod$, of $\num$.}

\begin{proof}
Define $\cod \from \Nat_+ \to \Nat_+$ as follows.
\display{10}{\begin{split} \cod(1) &= \min \{ k \mid 1 < k \mand
    \proj 1 \num(k) = \proj 1 \num(1) \mand\\  &\strut \qquad \qquad
    \max \{ 1, \proj 1 \num(1) \} < \proj 2 \num(k) \}\\
    \cod(n+1) &= \min \{ k \mid \cod(n) < k \mand \proj 1 \num(k) = \proj
    1 \num(n+1) \\ &\strut \qquad \qquad \mand \max \{ \proj 1
    \num(k), \proj 2 \cpair(n) \} < \proj 2 \num(k) \}
\end{split}}

By \rslt{3}, $\cod$ is computable if it is total. That
$\cod$ is total is proved by induction on the hypothesis that
$\cod(n)$ converges. To prove the basis step, let $\proj 1 \num(1)=i$
and note that the equivalence class of $i$ is infinite. Therefore, for
infinitely many $j>i$, $(i,j) \in E$. Since $\num$ enumerates $E$,
there are some $(i,j) \in E$ and $h > 1$ such that $(i,j) =
\num(h)$, so $\{ k \mid 1 < k \mand \proj 1 \num(1) = \proj 1 \num(k) <
\proj 2 \num(k) \}$ is non empty. $\cod(1)$ is defined to be the
least element of this set, so $\cod(1)$ converges.

The proof of the induction step is parallel. Suppose that
$\cod(n)$ converges. Let $\proj 1 \num(n+1) = i$. The equivalence
class of $i$ is infinite, so, for infinitely many $j > \max \{ \proj 1
\num(k), \proj 2 \cpair(n) \}$, $(i,j) \in E$. For some such $j$, and
for some $h > \cod(n)$, $(i,j) = \num(h)$. Thus $h \in \{ k \mid
\cod(n) < k \mand \proj 1 \num(k) = \proj 1 \num(n+1) \mand \proj 2
\num(k) > \max \{ \proj 1 \num(k), \proj 2 \cpair(n) \} \}$. Since
$\{ k \mid \proj 1 \num(k) = \proj 1 \num(n+1) \mand \proj 2 \num(k) >
\max \{ \proj 1 \num(k),\\ \proj 2 \cpair(n) \} \}$ is non empty, it
has a least element, so $\cod(n+1)$ converges. By the principle of
induction, then, $\cod$ is total, and therefore $\cod$ is computable.

Clearly \eqn{10} specifies a function that satisfies definition
\eqn{8} of a coding.
\end{proof}

\propo{8}{For a semi-decidable equivalence relation, $E$,
  the following three conditions are equivalent.
\be
\setcounter{enumi}{\value{equation}}
\meti{11}{Some enumeration of $E$ has a coding;}
\meti{12}{$E$ is \coarsep;}
\meti{13}{Every enumeration of $E$ has a coding.}
\setcounter{equation}{\value{enumi}}
\ee}

\begin{proof}
Condition \eqn{11} implies condition \eqn{12} by
\rslt{6}. Condition \eqn{12} implies condition \eqn{13} by
\rslt{7}. Since every semi-decidable equivalence relation has an
enumeration, condition \eqn{13} implies condition \eqn{11}.
\end{proof}

\Section{4}{Three relations defined from a coding}

Throughout sections \ref{sec:4}--\ref{sec:7}, it will be assumed that
\begin{enum}
\meti{14}{$\cod$ codes $\num$, an
enumeration of $E$, an \coarse equivalence relation.}
\end{enum}
Two new relations---$\rcod$ and $\scod$---will be defined in this
section, and they will be shown to be decidable. It will be shown that
$\rcod \subseteq E$, $\scod \subseteq E$, and $E = \rcod \scod \inv
\rcod$. A third relation, $\tcod$, will also be defined. It will be
shown that $\tcod$ is decidable, $\tcod \subseteq E$, and $\scod \cup
\tcod$ is symmetric and transitive, laying the groundwork for
extending the result that $E = \rcod \scod \inv \rcod$ to a result
that every semi-decidable, \coarse equivalence relation the lattice join of
decidable, \coarsep, equivalence relations.

Define
\display{15}{ \rcod = \cpair(\Nat_+)}
and
\display{16}{\scod = \{ (\proj 2 \cpair(m),
  \proj 2 \cpair(n)) \mid \num(n) = \inv{(\num(m))} \}}

\lema{9}{$\rcod$ and $\scod$ are decidable.} 

\begin{proof} First, consider $\rcod$. Since $n < \proj 2
  \cpair(n)$ (\rslt{5}),
\display{17}{\rcod = \{ (i,j) \mid \exists n \! < \! j \; (i,j) = \cpair(n)
  \}}
Since $\cpair$ is computable (\rslt{4}),
\display{18}{\{(i,j,k,n) \mid n < k \mand (i,j) = \cpair(n) \}}
 is a decidable relation. Decidable relations are closed under bounded
 quantification.\footnote{In the proofs of \rslt{15} and \rslt{23} below, it will be
   important that the bound may be the value of a computable function
   of variables of the relation. For example, if $R$ is decidable and
   $f$ is computable, then $S(x,y) \iff \exists z \! < \! f(x,y) \,
   R(x,y,z)$ defines a decidable relation, $S$.}
 (Cf.\ \citet[p.\ 311]{Rogers-1967}.) Therefore
\display{19}{\{ (i,j,k) \mid \exists n \! < \! k \; (i,j) = \cpair(n) \}}
 is decidable. Let $f$ be the characteristic function of $\{ (i,j,k) \mid
  \exists n \! < \! k \; (i,j) = \cpair(n) \}$, and let $g(i,j) =
  (i,j,j)$. Both $f$ and $g$ are computable, so $fg$, the characteristic
  function of $\rcod$, is computable. That is, $\rcod$ is decidable.

The decision procedure for $\scod$ is similar. $\scod$ is defined by
\eqn{16} which, setting $i = \proj 2 \cpair(m)$ and $j = \proj 2
\cpair(n)$, is equivalent to
\display{20}{\scod = \{ (i,j) \mid \exists m \, \exists
n \, [i = \proj 2 \cpair(m) ]\mand j = \proj 2 \cpair(n) \mand
 \num(n) = \inv{(\num(m))} \}}
By \rslt{5}, the existential quantifiers in \eqn{20} can be replaced by
bounded existential quantifiers, yielding
\display{21}{\begin{gathered}
\scod = \{ (i,j) \mid \exists m \! < \! i \; \exists
n \! < \!j \; [i = \proj 2 \cpair(m) \mand\\ j = \proj 2 \cpair(n) \mand
 \num(n) = \inv{(\num(m))}] \}
\end{gathered}}
The proof that $\scod$ is decidable from equation \eqn{21} is parallel to that
for $\rcod$ from equation \eqn{19}.
\end{proof}

\lema{10}{$\rcod \subseteq E$ and $\scod \subseteq E$.}

\begin{proof}
For $e \in \Nat \times \Nat$, $e \in \rcod$ if and only if, for some
$n$, $e = \cpair(n)$. Thus there is some number (specifically, $m =
\cod(n)$), such that $e = \num(m)$, so $e \in E$.

If $(i,j) \in \scod$, then, for some $m,\, n,\, p,\, q$, 
\display{22}{i = \proj 2 \cpair(m) \quad \num(m) = (p, q) \qquad
j = \proj 2 \cpair(n) \quad \num(n) = (q,p)}
There are $h$ and $k$ such that $(h,i) = \cpair(m)$ and $(k,j) = \cpair(n)$.
Since $\forall r \; \proj 1 \cpair(r) = \proj 1 \num(r)$, $h = p$ and
$k = q$. Therefore $(h,k) = (p,q) = \num(m) \in E$.

Since $\cpair$ enumerates $\rcod \subseteq E$, $(h,i) \in E$. Since
$E$ is symmetric, $(i,h) \in E$. Since also $(h,k) \in E$ and $E$ is
transitive, $(i,k) \in E$. Finally, since also $(k,j) \in E$, $(i,j)
\in E$.
\end{proof}

\propo{11}{$E = \rcod \scod \inv \rcod$. An \coarse equivalence
relation is semi-decidable if, and only if, it is a composition
of decidable relations.}

\begin{proof}
By \rslt{10}, $\rcod \scod \inv \rcod \subseteq E$. To show that $E
\subseteq \rcod \scod \inv \rcod$, suppose that $e = (i,j) \in
E$. Let $e = \num(m)$ and let $\inv e = \num(n)$. Let $(i,h) =
\cpair(m)$ and let $(j,k) = \cpair(n)$. Then $(i,h) \in
\rcod$, $(h,k) \in \scod$, and $(k,j) \in \inv \rcod$, so $e \in
\rcod \scod \inv \rcod$.

That every composition of decidable relations is semi-decidable, is a
routine result of computability theory.\footnote{Cf.\ \citet[problem
    5-18]{Rogers-1967}. Rogers uses `relative product' to denote the
  composition of two relations.} An \coarse, semi-decidable
equivalence relation is codable by \rslt{8}, so $E = \rcod \scod \inv
\rcod$ establishes that every \coarse, semi-decidable relation is a
composition of decidable relations.
\end{proof}

To prove that a semi-decidable, \coarse equivalence relation is the
lattice join of decidable equivalence relations, will require that
there should be a decidable, transitive, symmetric relation, $U$, such
that $\scod \subseteq U \subseteq E$. This relation will be obtained
by taking $U$ to be the union of $\scod$ with the following relation.
\display{23}{\tcod = \{ (i,j) \mid \exists m \! < \! i \; \exists
n \! < \!j \; [i = \proj 2 \cpair(m) \mand j = \proj 2 \cpair(n) \mand
 \num(n) = \num(m)] \}}
Also define the \emph{identity relation,}
\display{24}{\id = \{ (n,n) \mid n \in \Nat \}}

\lema{12}{$\tcod$ is decidable and $\tcod \subseteq E$. $\scod \cup
  \tcod$ is symmetric and transitive. Thus $\scod \cup \tcod \cup \id
  \subseteq E$ is an equivalence relation. $\scod \cup \tcod$ and
  $\scod \cup \tcod \cup \id$ are decidable relations.}

\begin{proof}
Since $\scod$ and $\tcod$ are subsets of $E$ (\rslt{10}) and also $\id
\subseteq E$, $\scod \cup \tcod \cup \id \subseteq E$.

The proof that $\tcod$ is decidable closely follows the corresponding
proof for $\scod$. To prove that $\tcod \subseteq E$, suppose that
$(i,j) \in \tcod$.Then, for some $m,\, n,\, p,\, q$, 
\display{25}{i = \proj 2 \cpair(m) \quad j = \proj 2 \cpair(n) \quad
 \num(m) = \num(n) = (p,q)}
There are $h$ and $k$ such that $(h,i) = \cpair(m)$ and $(k,j) = \cpair(n)$.
Since $\forall r \; \proj 1 \cpair(r) = \proj 1 \num(r)$, $h = k =
p$. That is, $(p,i) = \cpair(m)$ and $(p,j) = \cpair(n)$, so $(p,i)
\in E$ and $(p,j) \in E$. Since $E$ is symmetric and transitive,
$(i.j) \in E$.

Substituting $j$ and $n$ for $i$ and $m$ respectively in equations
\eqn{16} and \eqn{23} results in equivalent expressions.  Thus $\scod$
and $\tcod$, and consequently $\scod \cup \tcod$, are symmetric. To
see that $\scod \cup \tcod$ is transitive, suppose that $(i,j) \in
\scod \cup \tcod$ and $(j,k) \in \scod \cup \tcod$. If exactly one of
$(i,j)$ and $(j,k)$ is in $\scod$, then $(i,k) \in \scod$. Otherwise,
$(i,k) \in \tcod$.

Decidable relations are closed under union. Since $\scod$ is decidable
(\rslt{9}) and $\tcod$ is also decidable (by parallel reasoning) and
$\id$ is decidable (obvious), $\scod \cup \tcod$ and $\scod \cup \tcod
\cup \id$ are decidable. Since $\scod \cup \tcod$ is symmetric and
transitive, and taking their union with $\id$ preserves those
properties and makes the resulting relation reflexive, $\scod \cup
\tcod \cup \id$ is a decidable equivalence relation.
\end{proof}

\Section{5}{$\tc{(\rcod \cup \inv \rcod)}$ is decidable}

It will be shown that $\tc{(\rcod \cup \inv \rcod)}$ is a decidable,
\coarsep, equivalence relation. The first step is to show that it is
decidable. By \rslt{3}, $(i,j) \in \tc{(\rcod \cup \inv \rcod)}$ iff
there is a walk from $i$ to $j$ in $\rcod$.

\lema{13}{Suppose that $i \neq j$ and that $(x_1,\dotsc,x_n)$ is a walk
  of minimal length from $i$ to $j$ in $\rcod$, where $\cpair$
  enumerates $\rcod$. Then, for some $k$, $0 \le k < n$ and $\forall t
  \! \le k \; x_t < 0$ and $\forall t \! > k \; x_t > 0$.}

\begin{proof}
Suppose that $i \neq j$ and that $(x_1,\dotsc,x_n)$ is a walk of
minimal length from $i$ to $j$ in $\rcod$. If $\forall t \; x_t > 0$,
then $k = 0$ satisfies the required condition. If $\forall t \; x_t <
0$, then $k = n$ satisfies the condition. If $\exists p \, \exists q
\, [x_q < 0 \mand x_p > 0]$, then define $q^* = \max \{ q \mid x_q < 0
\}$. If $\forall q \! < \! q^* \; x_q < 0$, then $k = q^*$ satisfies
the condition. Otherwise, set $p^* = \max \{ p \mid p < q^* \mand x_p
> 0 \}$. By the definition of a walk, then, $\proj 2 \cpair(x_{p^*}) =
\proj 2 \cpair(-x_{p^*+1})$. Since $\proj 2 \cod$ is strictly
increasing, (\rslt{5}), $x_{p^*} = -x_{p^*+1}$. It follows that $n
>2$, since otherwise $i = \proj 1 \cpair(x_{p^*}) = \proj 1
\cpair(-x_{p^*+1}) = j$, contrary to hypothesis. Since $n > 2$,
deleting $x_{p^*}$ and $x_{p^*+1}$ from the walk does not delete the
entire walk, so $(x_1,\dotsc,x_{p^*-1},x_{p^*+2},\dotsc,x_n)$ is a
walk from $i$ to $j$ in $\rcod$. This contradicts the minimality of the
length of $(x_1,\dotsc,x_n)$. Thus, since $\exists p \, [p < q^* \mand
  x_p > 0]$ is impossible, $k = q^*$ satisfies the required condition.
\end{proof}

Computable functions can be defined that respectively associate numbers with
integer sequences of arbitrary positive length and bound the smallest
representing number of a sequence. Specifically,

\lema{14}{There exist computable functions $\seq \from \Nat_+ \times
  \Nat \to \Z$ and $\bound \from \Nat_+ \to \Nat_+$ such that, for
  every $n > 0$, $(x_1,\dotsc,x_n) \in (\Z \setminus \{ 0 \})^n$, $k
  \ge n$, and $k \ge \max \{ \lvert x_i \rvert \mid 1 \le i \le n \}$,
  \display{26}{ \exists z \le \bound (k) \enspace [ \seq(z,0) = n
      \mand \forall i \! < \!  n \; \seq(z,i+1) = x_i ] } }

\begin{proof}
  Begin by defining an injection, $f \from \bigcup_{n > 0} (\Z
  \setminus \{ 0 \})^n \to \Nat$, as follows. Represent $x$, a
  non-zero integer, by the string consisting of $\vert x \vert$
  ocurrences of `0', preceded by `1' if $x < 0$ and `11' if $x >
  0$. Represent a sequence of non-zero integers by the concatenation
  of their strings. If `$d^1\cdots d^n$' (a concatenation of $n$
  binary-digit strings, $d^i$, representing the respective $x_i$) is
  the binary-digit string representing $(x_1,\dotsc,x_n) \in (\Z
  \setminus \{ 0 \})^n$, and if `$d^1\cdots d^n$' is
  $\text{`}b_1\text{'}\cdots \text{`}b_t\text{'}$ (where each
  $\text{`}b_k\text{'}$ is `0' or `1') then define $f(x_1,\dotsc,x_n)
  = \sum_{s < t} b_{n-s} \, 2^s$.

Define $\seq \from \Nat_+ \times \Nat \to \Z$ by
\display{27}{\seq(z,i) = \begin{cases}
0 & \text{if $z$ is not in the range of $f$, and}\\
& \text{\strut \quad if $z = f(x_1,\dotsc,x_n)$, then}\\
n & \text{if\ } i = 0;\\
x_i & \text{if\ } 0 < i \le n;\\
0 & \text{if\ } i>n
\end{cases}
}
This definition works because the set of numbers corresponding to
nonzero-integer sequences, and also the relation, ``$x$ is the
$i^{\text{th}}$ element of the nonzero-integer sequence corresponding
to $z$'', are decidable. That is true, in turn, because occurrences of
`10' and `11' occur exactly at the beginnings substrings of a
concatenated string that represent the respective $x_i$ and because
the range of $f$ is decidable.\footnote{A string of digits represents
  a sequence of non-zero integers iff it satisfies three conditions:
  the first digit must be `1', the last digit must be `0', and there must
  be no substring of form `111'. These conditions can be expressed by
  bounded-quantifier formulae in base-2 arithmetic, so the range of
  $f$ is decidable. A detailed proof of this decidability assertion
  would proceed along the general lines of \citet[definition V.3.5 and
    proposition V.3.30]{HajekPudlak-1998}.} The $z$ satisfying
\eqn{26} can be bounded because, if the binary-digit string,
$\text{`}b_1\text{'}\cdots \text{`}b_t\text{'}$, encodes $(x_1,
\dotsc, x_n)$, then $t \le \sum_{i=1}^n (2 + x_i)$, so
$f(x_1,\dotsc,x_n) \le \sum_{s=0}^{\sum_{i=1}^n (2 + x_i)} 2^s < 2^{1
  + n (2+ \max \{ \lvert x_i \rvert \mid 1 \le i \le n \})}$. That is,
if $\bound(k) = 2^{1 + k (2 + k)}$ and $n \le k$ and $\max \{ \lvert
x_i \rvert \mid 1 \le i \le n \} \le k$, then $f(x_1,\dotsc,x_n) <
\bound(k)$ and $z = f(x_1,\dotsc,x_n)$ witnesses \eqn{26}.
\end{proof}

\lema{15}{$\tc{(\rcod \cup \inv \rcod)}$ is a decidable relation}

\begin{proof}
$\rcod$ is decidable by \rslt{9}, so $(\rcod \cup \inv \rcod)$ is
  decidable. By \rslt{3}, $(i,j) \in \tc{(\rcod \cup \inv \rcod)}$ iff
  there is a walk from $i$ to $j$ in $\rcod$. Let $(x_1, \dotsc, x_n)$
  be such a walk of minimal length. By \rslt{13}, for some $k$, $0 \le
  k < n$ and $\forall i \! \le k \; x_i < 0$ and $\forall j \! > \! k
  \; x_j > 0$. This means that $i = \proj 2 \cpair(-x_1)$; that, for
  $h \le k$, $\proj 1 \cpair(-x_h) = \proj 2 \cpair(-x_{h+1})$; that
  $\proj 1 \cpair(x_{k+1}) = \proj 2 \cpair(-x_k)$; that, for $h \ge
  k$, $\proj 2 \cpair(x_h) = \proj 1 \cpair(x_{h+1})$; and that $\proj
  2 \cpair(x_n) = j$. These facts, together with the fact that $\proj
  1 \cpair( \vert x_h \vert) < \proj 2 \cpair( \vert x_h \vert)$
  (definition \eqn{8}), have two implications. First, $\proj 2
  \cpair(-x_h) \le i$ for $h \le k$ and $\proj 2 \cpair(-x_h) \le j$
  for $h > k$. Thus, for all $h$, $\proj 2 \cpair( \vert x_h \vert )
  \le \max \{ i,j \} \le i+j$. Second, $n \le i+j$. Thus, in view of
  \rslt{14}, it is clear that the set of $(i,j)$ such that there is a
  walk from $i$ to $j$ in $\rcod$ is defined from $\rcod$ by the
  bounded-quantifier formula
  \display{28}{\begin{gathered}
      \exists x \! \le \! \bound(i+j) \, [\tail(\seq(x,1)) = i \mand
        \head(\seq(x,\seq(x,0))) = j \mand\\ \forall k \! < \!
        \seq(x,0) - 1 \, [ \head(\seq(x, k+1)) = \tail(\seq(x, k+2)
          \mand\\ (\seq(x, k+1), \seq(x, k+2)) \in \rcod \cup \inv
          \rcod ]]
  \end{gathered}}
  Thus that set---which is $\tc{(\rcod \cup \inv \rcod)}$ by
  \rslt{3}---is decidable.
\end{proof}

\Section{6}{Decidable and semi-decidable \coarse relations}

The main result of this section, \rslt{17}, will provide a partial
answer to the question, under what conditions is a semi-decidable
relation the lattice join of of decidable relations?

\lema{16}{$\tc{(\rcod \cup \inv \rcod)}$ is a decidable, \coarse
  equivalence relation.}

\begin{proof}
Clearly $\rcod \cup \inv \rcod$ is symmetric, and the transitive
closure of a symmetric relation is also symmetric as well as being
transitive. For every $i$, there is some $n$ such that $(i,i) =
\num(n)$, so $\cpair(n) \in \rcod$, $\inv{(\cpair(n))} \in \inv
\rcod$, $i = \proj 1 \cpair(n)$ and therefore $(i,i) \in \rcod \inv
\rcod \subseteq \tc{(\rcod \cup \inv \rcod)}$. That is, $\tc{(\rcod
  \cup \inv \rcod)}$ is an equivalence relation.  $\tc{(\rcod \cup
  \inv \rcod)}$ is $\coarsep$ by \rslt{6}, and is decidable by
\rslt{15}.
\end{proof}

Recall that, if $\eqva$ and $\eqvb$ are equivalence relations on $\Nat$, then
their join in the lattice of equivalence relations is defined by
\display{29}{\eqva \vee \eqvb = \tc{(\eqva \cup \eqvb)}}

\propo{17}{If $E$ is a semi-decidable, \coarse equivalence relation,
  then there are decidable equivalence relations,
  $\eqva$ and $\eqvb$, such that $\eqva$ is \coarse and
\display{30}{E = \eqva\eqvb\eqva = \eqva \vee \eqvb}
Specifically,
\display{31}{\eqva = \tc{(\rcod \cup \inv \rcod)} \qquad \eqvb = \scod \cup
\tcod \cup \id}
}

\begin{proof}
By definition \eqn{31} and \rslt{16}, $\eqva$ is a decidable, \coarse
equivalence relation. By \rslt{12}, $\eqvb$ is a decidable
equivalence relation. They are both sub-relations of $E$ so $\eqva
\vee \eqvb \subseteq E$. By \rslt{11}, $E = \rcod \scod \inv \rcod
\subseteq \eqva\eqvb\eqva \subseteq \eqva \vee \eqvb \subseteq E$,
which is equivalent to \eqn{30}.
\end{proof}

\Section{7}{Decidable \coarse equivalence relations generate\\ 
the semi-decidable \coarse equivalence relations}

It is of interest to strengthen \rslt{17} to assert that $\eqvb$, as
well as $\eqva$, is \coarsep. That is,

\propo{18}{If $E$ is a semi-decidable, \coarse equivalence relation,
  then there are decidable, \coarse equivalence relations,
  $\eqva$ and $\eqvb$, such that condition \eqn{30} holds.}

\Rslt{18} implies that the set of decidable, \coarse equivalence relations is
large enough to generate the semi-decidable \coarse equivalence
relations as a subset of the semigroup of binary relations, and as an
upper semilattice.

If $\eqvb$ is defined as in \eqn{31}, then $[i]_\eqvb$ is not
guaranteed to be infinite for $i \in \proj 2 \cpair(\Nat_+)$, and
$[i]_\eqvb = \{ i \}$ for $i \notin \proj 2
\cpair(\Nat_+)$.\footnote{Since $\proj 2 \cpair$ is strictly
  increasing, and since, because $E$ is \coarsep, there are infinitely
  many $n$ such that $\proj 2 \num(n+1) = \proj 2 \cpair(1)$, $\Nat_+
  \setminus \proj 2 \cpair(\Nat_+)$ is infinite.} There does not seem
to be any ad hoc adjustment of the definition of $\eqvb$ in \eqn{31}
that will yield, for arbitrary $\cod$, a decidable, \coarse
equivalence relation that includes $\scod \cup \tcod$, which will be
required required in order to transform the proof of \rslt{17} into a
proof of \rslt{18}.

\Rslt{17} and the results on which it depends have been proved by
starting with an arbitrary coding, $\cod$, of an arbitrary enumeration
$\num$, of $E$. In order to obtain the decidable, \coarse equivalence
relation that is required to prove \rslt{18}, $\num$ and $\cod$ will be
constructed, each of which has a specific property. For every $e \in
E$, $\inv \num(e)$ will be infinite. For each $i \notin \proj 2
\cpair(\Nat_+)$, $[i]_E \setminus \proj 2 \cod(\Nat_+)$ will be
infinite. These two properties will ensure that $\eqva$ is \coarse and
enable a decidable, \coarsep, symmetric and transitive relation,
$\eqvc$ to be constructed that will play a cognate role to that of
$\id$ in definition \eqn{31} of $\eqvb$.

\lema{19}{There is an enumeration, $\num$, of $E$ (a semi-decidable, coarse
  equivalence relation) such that
\display{32}{\text{for all\ } (i,j) \in E,\; \inv \num(i,j) \text{\ is infinite.}}
}

\begin{proof}
The mapping $(m,n) \mapsto 2^m(2n+1)$ is a bijection of $\Nat \times
\Nat$ to $\Nat_+$. Since $E$ is semi-decidable, it has some
enumeration, $\numtwo$. An enumeration satisfying \eqn{32} is defined by
$\num(2^m(2n+1)) = \numtwo(n+1)$.
\end{proof}

\lema{20}{There are an enumeration, $\num$,  of $E$ (a semi-decidable,
  coarse equivalence relation), and codings, $\cod$ and
  $\codtwo$, of $\num$, such that $\num$ satisfies \eqn{32} and
\display{33}{\text{for all\ } n \in \Nat_+,\; \cod(n) < \codtwo(n) < \cod(n+1)}
}

\begin{proof}
  The enumeration exists by \rslt{19}. Intuitively, the codings are
  constructed by interleaving two, concurrent recursions, each
resembling the one used to prove \rslt{7}. At each stage, $n$, implicitly first
$\cod(n)$ is constructed and then $\codtwo(n)$ is constructed and then
the results are merged by use of a pairing function, $(i,j) \mapsto
2^i(2j+1)$. That is, define
\display{34}{\merge(n) = 2^{\cod(n)}(2\codtwo(n)+1)}

To describe this procedure within a single, recursive, definition of
$\merge$, define
\display{35}{\begin{gathered}
    \Proj 1(k) = \max \{ m \mid \; 2^m \vert k \} \qquad \Proj 2(k) =
    [(k/\Proj 1(k))-1]/2\\
    \cod = \Proj 1 \merge \qquad \codtwo = \Proj 2 \merge \qquad
    \cpairtwo = \num \codtwo
\end{gathered}}
Now, $\merge \from \Nat_+ \to \Nat_+$ will be defined so that $\cod$
and $\codtwo$, derived from $\merge$ via \eqn{35}, will each satisfy the
definition \eqn{10} of a coding, and so that they will be interleaved as
specified in \eqn{33}.
\display{36}{\begin{split} \merge(1) &= \min \{ k \mid 1 < \min \{ \Proj
    1(k), \Proj 2(k) \} \mand\\
 &\proj 1 \num(1) = \proj 1 \num(\Proj 1(k)) = \proj 1 \num(\Proj 2(k)) \mand\\
 &\max \{ 1, \proj 1 \num(1) \} <
    \proj 2 \num(\Proj 1(k)) <  \proj 2 \num(\Proj 1(k)) \}\\
\merge(n+1) = \min &\{ k \mid \codtwo(n) < \Proj 1(k) < \Proj 2(k) \mand\\
    \proj 1 \num(n+1) &= \proj 1 \num(\Proj 1(k)) = \proj 1 \num(\Proj 2(k)) \mand\\
    \max \{ \proj 1 \num(\Proj 1(k)), &\proj 1 \num(\Proj 2(k)), \proj
    2 \cpairtwo(n) \} < \proj 2 \num(\Proj 1(k)) < \proj 2 \num(\Proj 2(k)) \}
\end{split}}
By a parallel argument to the proof of \rslt{7}, the fact that $E$ is
\coarse implies that \eqn{36} defines a total function. Clearly $\cod$
and $\codtwo$, defined by \eqn{35} and \eqn{36}, each satisfy definition
\eqn{8} of a coding, and they are related according to \eqn{33}.
\end{proof}

\lema{21}{If $\cod$ codes $\num$, an enumeration of $E$ that satisfies
  \eqn{32}, then $\scod \cup \tcod$ is an \coarsep, symmetric and
  transitive relation.}

\begin{proof}
By \rslt{12}, $\scod \cup \tcod$ is symmetric and transitive. Let $i
\in \field(\scod \cup \tcod)$. It must be shown that $[i]_{\scod \cup
  \tcod}$ is infinite.

For some $n$, $i = \proj 2 \cpair(n)$. Let $\cpair(n) = (j,i)$. Let $M
= \inv \num(j,i)$. By \eqn{32}, $M$ is infinite. Since $\proj 2
\cpair$ is strictly increasing, $\proj 2 \cpair(M)$ is infinite. For
each $m \in M$, $(\proj 2 \cpair(m), \proj 2 \cpair(n)) \in \tcod$, so
$\proj 2 \cpair(M) \subseteq [i]_{\tcod} \subseteq [i]_{\scod \cup
  \tcod}$, so $[i]_{\scod \cup \tcod}$ is infinite.
\end{proof}

The next step is to include $\scod \cup \tcod$ in a decidable,
\coarsep, equivalence relation on $\Nat$. This will be done by
appealing to \rslt{1}. Setting $H = \scod \cup \tcod$, it is
sufficient to find another relation that satisfies the conditions
specified in the following lemma.

\lema{22}{There is a sub-relation of $E$ that is decidable,
\coarsep, symmetric and transitive, and has
$\Nat \setminus \field(\scod \cup \tcod)$ as its field.}

\begin{proof}
The following relation will be proved to satisfy the specified
conditions.
\display{37}{J = ((\Nat \setminus \field(\scod \cup \tcod)) \times (\Nat
  \setminus \field(\scod \cup \tcod)) \cap \tc{(\rcod[\codtwo] \cup \inv
            {\rcod[\codtwo]})}}

$\field(\scod \cup \tcod)$ is Turing reducible to $\scod \cup \tcod$
because $i \in \field(\scod \cup \tcod) \iff (i,i) \in \scod \cup
\tcod$ (\rslt{1}).  $\rcod[\codtwo] \cup \inv {\rcod[\codtwo]}$ is
decidable (\rslt{12}), so $\field(\scod \cup \tcod)$ is
decidable. Therefore $(\Nat \setminus \field(\scod \cup \tcod)) \times
(\Nat \setminus \field(\scod \cup \tcod))$ is decidable,
also. $\tc{(\rcod[\codtwo] \cup \inv {\rcod[\codtwo]})}$ is decidable
(\rslt{15}), so $J$ is decidable. $J \subseteq E$ because
$\tc{(\rcod[\codtwo] \cup \inv {\rcod[\codtwo]})} \subseteq E$. Being
an intersection of symmetric, transitive relations, $J$ shares those
properties.

It remains to be proved that $J$ is \coarse and that $J$ has $\Nat
\setminus \field(\scod \cup \tcod)$ as its field. A single argument
establishes both of these facts. It follows directly from \eqn{37}
that $\field(J) \cap \field(\scod \cup \tcod) = \emptyset$, so, to
prove that $\field(J) = \Nat \setminus \field(\scod \cup \tcod)$, it
is sufficient to prove that $\Nat \setminus \field(\scod \cup \tcod)
\subseteq \field(J)$. To that end, suppose that $i \notin \field(\scod
\cup \tcod)$, and suppose that $(i,i) = \num(k)$. By \rslt{6} and
\eqn{15}, $\{ \cpairtwo\codtwo^n(k) \mid n \in \Nat \} = \{ i \} \times
\{ \proj 2 \cpairtwo\codtwo^n(k) \mid n \in \Nat \}$ is an infinite
subset of $\rcod[\codtwo]$, hence of $\tc{(\rcod[\codtwo] \cup \inv
  {\rcod[\codtwo]})}$.  \Rslt{6} and condition \eqn{33} imply that $\{
\proj 2 \cpairtwo\codtwo^n(k) \mid n \in \Nat \} \subseteq \Nat \setminus
\field(\scod \cup \tcod))$, so $\{ i \} \times \{ \proj 2
\cpairtwo\codtwo^n(k) \mid n \in \Nat \} \subseteq (\Nat \setminus
\field(\scod \cup \tcod)) \times (\Nat \setminus \field(\scod \cup
\tcod)$. Thus $\{ i \} \times \{ \proj 2 \cpairtwo\codtwo^n(k) \mid n \in
\Nat \} \subseteq J$ and $[i]_J$ is infinite. Since $[i]_J \neq
\emptyset$, $i \in \field(J)$.
\end{proof}

The main result of this article, \rslt{18}, follows directly from preceding lemmas. 

\begin{proof}[Proof of \rslt{18}]
Let  $\num$, $\cod$, and $\codtwo$ be functions such as are guaranteed by
\rslt{20} to exist. Define $J$ according to \eqn{37}, and also define
\display{38}{\eqva = \tc{(\rcod \cup \inv \rcod)} \qquad H = \scod \cup \tcod
\qquad \eqvb = H \cup J}
Lemmas \rsltref{1}, \rsltref{15}, \rsltref{21}, and \rsltref{22}
establish that $\eqva \subseteq E$ and $\eqvb \subseteq E$ are
decidable, \coarsep, equivalence relations. That $E = \eqva\eqvb\eqva
= \eqva \vee \eqvb$ follows as in the proof of \rslt{17}.
\end{proof}

\Section{8}{semi-decidable \fine relations}

Two examples of \fine relations that are semi-decidable but not decidable will
now be constructed. One of the examples has a representation of form
$E = \eqva \vee \eqvb$, where $\eqva$ and $\eqvb$ are decidable, while the other
cannot be so represented.

\propo{23}{There exist decidable equivalence relations, $\eqva$ and $\eqvb$,
  such that $\eqva \vee \eqvb$ is a semi-decidable, \fine equivalence
  relation that is not decidable, and $\eqva \vee \eqvb = \eqva\eqvb\eqva$.}

\begin{proof}
Let $\numtwo \from \Nat_+ \to A$ be an enumeration of $A$, a
semi-decidable, but not decidable, subset of $\Nat$. Assume, in
accordance with \citet[problem 5-2, p.~73]{Rogers-1967}, that
$\numtwo$ enumerates $A$ without repetitions. Define two
equivalence relations, $\eqva$ and $\eqvb$, by
\display{39}{(i,j) \in \eqva \iff \exists n \! \le \! \max \{ i,j \}
  \, [i=j \text{\ or\ } \{ i,j \} = \{ 2n, 3^{\numtwo(n)} \} ]}
and
\display{40}{(i,j) \in \eqvb \iff \exists n \! \le \! \max \{ i,j \}
  \, [i=j \text{\ or\ } \{ i,j \} = \{ 2n, 5^{\numtwo(n)} \} ]}
Since $\eqva$ and $\eqvb$ are defined from equality and
computable functions by bounded-quantifier sentences, they are
decidable. They are symmetric and reflexive and, since $\numtwo$ is an
enumeration without repetition, each equivalence class is either a
singleton or a pair. Thus, by \rslt{2}, $\eqva$ and $\eqvb$ are
equivalence relations.

Since $\numtwo$ is bijective, the partition corresponding to $\eqva$
is $\{ \{ i,j \} \mid [i=j \mand \neg \exists n$ $[i = 2n
    \text{\ or\ } i = 3^{\numtwo(n)} ]] \text{\ or\ }  \exists n \,
  \{ i,j \} = \{ 2n, 3^{\numtwo(n)} \} \}$ and the partition corresponding to
$\eqvb$ is $\{ \{ i,j \} \mid [i=j \mand \neg \exists n [i = 2n
    \text{\ or\ } i = 5^{\numtwo(n)} ]] \text{\ or\ }  \exists n \,
  \{ i,j \} = \{ 2n, 5^{\numtwo(n)} \} \}$

Suppose that $(i,j) \in \eqva\eqvb\eqva$. Then there are $h$ and $k$
such that $(i,h) \in \eqva$, $(h,k) \in \eqvb$, and $(k,j) \in
\eqva$. If $i = h$, then $(i,j) \in \eqvb\eqva$. If $i \neq j$, then
either $h = k$ or $h \neq k$. If $h = k$, then $(i,j) \in
\eqva\id\eqva = \eqva \subseteq \eqva\eqvb$. If $h \neq k$, then, for
some $n$, $i = 3^{\numtwo(n)}$ and $h = 2n$ and $k =
5^{\numtwo(n)}$. In that case, since $(k,j) \in \eqva$, $j = k$, so
$(i,j) \in \eqva\eqvb$. Thus $\eqva\eqvb\eqva \subseteq \eqva\eqvb
\cup \eqvb\eqva = \eqva\eqvb\id \cup \id\eqvb\eqva \subseteq
\eqva\eqvb\eqva \cup \eqva\eqvb\eqva = \eqva\eqvb\eqva$. A parallel
argument shows that $\eqvb\eqva\eqvb \subseteq \eqva\eqvb \cup
\eqvb\eqvb \subseteq \eqva\eqvb\eqva$.  Therefore $\eqva\eqvb\eqva
\cup \eqvb\eqva\eqvb \subseteq \eqva\eqvb \cup \eqvb\eqva \subseteq
\eqva\eqvb\eqva$. That is,
\display{41}{\eqva\eqvb\eqva \cup \eqvb\eqva\eqvb = \eqva\eqvb \cup
  \eqvb\eqva = \eqva\eqvb\eqva}

Using the first identity in \eqn{41} to carry out the induction step, it
is seen by induction that $\forall n \, [n \ge 2 \implies
  \relpow{(\eqva \cup \eqvb)}{n} = \eqva \eqvb \cup \eqvb
  \eqva]$. Then, by equations \eqn{4}  and \eqn{29} and the second identity
in \eqn{41}, $\eqva \vee \eqvb = \eqva\eqvb\eqva$. Clearly the partition
corresponding to $\eqva \vee \eqvb$ is $\{ \{ i,j, k\} \mid [i=j=k
  \mand \neg \exists n \, [i = 2n \text{\ or\ } i = 3^{\numtwo(n)}
    \text{\ or\ } i = 5^{\numtwo(n)} ]] \text{\ or\ } 
  \exists n \, \{ i,j,k \} = \{ 2n, 3^{\numtwo(n)},
    5^{\numtwo(n)} \} \}$, so $\eqva \vee \eqvb$ is \finep.

Clearly $\eqva \vee \eqvb$ is semi-decidable. But
$A$ is Turing reducible to $\eqva \vee \eqvb$ by the equivalence $n
\in A \iff (3^n,5^n) \in \eqva \vee \eqvb$. Therefore $\eqva \vee
\eqvb$ not decidable.
\end{proof}

\propo{24}{There is a semi-decidable, \fine equivalence relation,
  $E$, that is not the transitive closure of a decidable relation.}

\begin{proof}
Let $A$ be a semi-decidable subset of $\Nat$ that is not decidable. 
Consider the semi-decidable equivalence relation, $E$, defined by
\display{42}{\begin{gathered}
(x,y) \in E \iff [x = y \text{\ or\ } [\min\{x,y\} \in 2\Nat\\
      \mand \min\{x,y\}/2 \in A \mand (y-x)^2 = 1]]
\end{gathered} }
$E$ corresponds to the partition, the equivalence classes of
which are defined by
\display{43}{[i] = \begin{cases}
\{ i \} & \text{if\ }  2n \le i \le 2n+1 \mand n \notin A;\\
\{ 2n, 2n+1 \} & \text{if\ }  2n \le i \le 2n+1 \mand n \in A.
\end{cases} }
$A$ is Turing reducible to $E$ because $n \in A \iff \{ 2n, 2n+1 \}\in
E$, so $E$ is not decidable.

A contradiction will be derived from the assumption that, for some
decidable $R$, $E = \tc R$. If so, then, since $\tc{(R \cup I)} = R
\cup I$ by \rslt{2}, $E = \tc R \subseteq \tc{(R \cup I)} = R \cup I
\subseteq E$. That is, $E = R \cup I$, which is impossible since $R
\cup I$ is decidable but $E$ is not decidable.
\end{proof}

\corol{25}{There is a semi-decidable, \fine equivalence relation, $E$,
  such that (a) there do not exist $n>1$ and decidable equivalence relations,
  $\eqvc_1, \dotsc,\eqvc_n$, that satisfy $E = \bigvee_{i \le n} \eqvc_i$, and (b)
  there do not exist $n>1$ and decidable equivalence relations,
  $\eqvc_1,\dotsc,\eqvc_n$ (not necessarily distinct), such that $E =
  \eqvc_1\cdots \eqvc_n$.}

\begin{proof}
Let $E$ be defined by \eqn{42}, with $A$ being a non decidable,
semi-decidable subset of $\Nat$. To prove (a), suppose that each
$\eqvc_i$, for $1 \le i \le n$< is a decidable equivalence
relation. Then $\bigcup_{i \le n} \eqvc_i$ is decidable and
$\bigvee_{i \le n} \eqvc_i = \tc{(\bigcup_{i \le n} \eqvc_i)}$.
By \rslt{24}, then, $E \neq \bigvee_{i \le n} \eqvc_i$. To prove (b)
  by contradiction, suppose that $\eqvc_1,\dotsc,\eqvc_n$ are
  decidable equivalence relations such that $E = \eqvc_1\cdots
  \eqvc_n$.  Because $\id \subseteq \eqvc_j$ for every $j$, $\eqvc_i
  \subseteq \eqvc_1\cdots \eqvc_n = E$ for every $i$. Thus $\bigcup_{i
    < n} \eqvc_i \subseteq E$. Also $\bigcup_{i < n} \eqvc_i$ is
  reflexive and symmetric, so, by \rslt{2}, $\tc{(\bigcup_{i < n}
    \eqvc_i)} = \bigcup_{i < n} \eqvc_i$. $\eqvc_1\cdots \eqvc_n
  \subseteq \relpow{(\bigcup_{i < n} \eqvc_i)}{n} \subseteq
  \tc{(\bigcup_{i < n} \eqvc_i)}$, so $E \subseteq \bigcup_{i < n}
  \eqvc_i \subseteq E$. That is, $E = \bigcup_{i < n} \eqvc_i$,
  contradicting $\bigcup_{i < n} \eqvc_i$ being decidable but $E$ not
  being so.
\end{proof}

\bibliographystyle{plainnat}

\begin{thebibliography}{5}
\providecommand{\natexlab}[1]{#1}
\providecommand{\url}[1]{\texttt{#1}}
\expandafter\ifx\csname urlstyle\endcsname\relax
  \providecommand{\doi}[1]{doi: #1}\else
  \providecommand{\doi}{doi: \begingroup \urlstyle{rm}\Url}\fi

\bibitem[Andrews et~al.(2017)Andrews, Badaev, and Sorbi]{AndrewsEtAl-2017}
U.~Andrews, S.~Badaev, and A.~Sorbi.
\newblock A survey on universal computably enumerable equivalence relations.
\newblock In A.~Day, M.~Fellows, N.~Greenberg, B.~Khoussainov, A.~Melnikov, and
  F.~Rosamond, editors, \emph{Computability and Complexity}, pages 418--451.
  Springer, 2017.

\bibitem[Gao and Gerdes(2001)]{GaoGerdes-2001}
Su~Gao and Peter Gerdes.
\newblock Computably enumerable equivalence relations.
\newblock \emph{Studia Logica}, 67\penalty0 (1):\penalty0 27--59, 2001.

\bibitem[H\'{a}jek and Pudlak(1998)]{HajekPudlak-1998}
P.~H\'{a}jek and P.~Pudlak.
\newblock \emph{Metamathematics of First-Order Arithmetic}.
\newblock Springer, 1998.

\bibitem[Kleene(1952)]{Kleene-1952}
Stephen~C. Kleene.
\newblock \emph{Introduction to Metamathematics}.
\newblock D.~van Nostrand Company, 1952.

\bibitem[Rogers~Jr.(1967)]{Rogers-1967}
Hartley Rogers~Jr.
\newblock \emph{The theory of recursive functions and effective computability}.
\newblock McGraw-Hill, 1967.

\end{thebibliography}

\end{document}